\documentclass[graybox]{svmult}

\usepackage{mathptmx}       
\usepackage{helvet}         
\usepackage{courier}        
\usepackage{type1cm}        
\usepackage{makeidx}         
\usepackage{graphicx}        
\usepackage{multicol,booktabs}        
\usepackage[bottom]{footmisc}
\usepackage{showlabels}
\usepackage{hyperref}
\usepackage{cite}
\usepackage{fancyvrb}        

\usepackage{amsmath}
\usepackage{amssymb}

\newcommand{\R}{\mathbb{R}}

\makeindex             

\begin{document}

\title*{HJB-POD feedback control for Navier-Stokes equations.}
\author{Alessandro Alla and Michael Hinze}
\institute{Alessandro Alla, Michael Hinze \at Department of Mathematics, Universit\"at Hamburg, Bundesstr. 55, 20146 Hamburg, Germany,\\ \email{ alessandro.alla/ michael.hinze @uni-hamburg.de} }
%
%
\maketitle

\abstract{In this report we present the approximation of an infinite horizon optimal control problem for the evolutive Navier-Stokes system. The method is based on a model reduction technique, using a POD approximation, coupled with a Hamilton-Jacobi equation which characterizes the value function of the corresponding control problem for the reduced system. Although the approximation schemes available for the  HJB are shown to be convergent for any dimension, in practice we need to restrict the dimension to rather small numbers and this limitation affects the accuracy of the POD approximation. We will present numerical tests for the control of the time-dependent Navier-Stokes system in two-dimensional spatial domains to illustrate our approach and to show the
effectiveness of the method.}
\\


\section{Introduction}
In this report we investigate an infinite horizon optimal control problem for the time-dependent Navier-Stokes equations (NSE). The basic ingredient of the method is the coupling between a proper orthogonal decomposition (POD) approximation of the NSE and a Dynamic Programming scheme for the stationary Hamilton-Jacobi equation characterizing the value function of the optimal control problem. Due to the curse of dimensionality, we need to restrict the dimension of the POD system to a rather small number (typically 4). This limitation naturally affects the accuracy of the POD approximation (see \cite{V07}), and, as a consequence, the problem class which we can treat with this technique. It is well known that the solution of the HJB equation is not an easy task from the numerical point of view since viscosity solutions of the HJB equation are usually just Lipschitz-continuous.  Optimal  control problems for ODEs are solved by Dynamic Programming (DP), both analytically and numerically (see \cite{BCD:97} for a general presentation of this theory). From the numerical point of view, this approach has been developed for many classical control problems obtaining convergence results and a-priori error estimates (see the recent book from Falcone and Ferretti \cite{FF12}).
We should mention that a first tentative approach to couple POD and HJB equations is proposed by  Atwell and King \cite{AK01} for the control of the 1D heat equation. Kunisch and Volkwein in \cite{KV99, KV01} extend this approach to diffusion dominated equations and, in particular, Kunisch, Volkwein and Xie in \cite{KVX04} apply HJB-POD feedback control to the viscous Burgers equation.  We also mention an adaptive POD technique for 1D advection dominated problems proposed by the first author and Falcone in \cite{AF12, AF13b}.\\
\noindent
The novelty in this paper consists in the control of the 2D nonlinear time dependent Navier-Stokes system by means of DP equations and the reduction of the nonlinear term with the {\em Discrete Empirical Interpolation Method} due to Chaturantabut and Sorensen in \cite{CS09}.\\
\noindent
The paper is organized as follows.We first present the optimal control problem in Section 2, then we describe the DP equation in Section 3. Proper orthogonal decomposition is summarized in Section 4 and, finally, the numearical tests are presented in Section 5.

\section{The optimal control problem}

In this section we describe the optimal control problem.  The gouverning equations are the two non-stationary dimensional unsteady Navier-Stokes equations. The flow in the bounded domain $\Omega\subset\R^2$ is characterized by the velocity field $y:\Omega\times[0, T]\rightarrow \R^2$ and by the pressure $p:\Omega\times[0, T]\rightarrow \R$. The Navier-Stokes equations are given by

\begin{equation}\label{NS}
\left.
\begin{aligned}
y_t-\nu\Delta y+(y\cdot  \nabla y)+\nabla p=\sum_{i=1}^Nb_i(x)u_i(t)&\quad&\mbox{ in } \Omega\times(0,T],\\
\nabla \cdot y=0&\quad&\mbox{ in } \Omega\times (0,T],\\
y(\cdot,0)=y_0&\quad&\mbox{ in } \Omega,\\
y(\cdot,t)=y_b&\quad&\mbox{ in } \partial\Omega\times (0,T),
\end{aligned}
\right\}
\end{equation}


\noindent
where the viscosity of the flow is given by the parameter $\nu>0$. The control signals are elements of $\mathcal{U}\equiv\{ u:[0,T]\rightarrow U,\, u(\cdot)\in L^\infty(0,T)\}$, where $U$ is a compact subset of $\mathbb{R}^m.$ Later we take $U$ as a discrete set. The initial value and the boundary values are denoted by $y_0$ and $y_b$, respectively. Finally, the functions $b_i(x):\Omega\rightarrow\R^2$ play the role of the so called shape functions.\\
\noindent
The cost functional we want to minimize is given by
\begin{equation}
J(u):=\int_0^\infty \left(\|y(\cdot,t;u)-\bar{y}\|^2_{L^2(\Omega)}+\alpha |u(t)|^2\right)e^{-\lambda t}\,dt,
\end{equation}
where $\bar{y}$ is the desired state which we choose as the mean flow, $\alpha\in\R^+$ and $\lambda>0$ is the discount factor. 
The optimal control problem, then, can be formulated as
\begin{equation}\label{ocp}
\min_{u\in\mathcal{U}} J(u)\mbox{ s. t. } y(u) \mbox{ satisfies } \eqref{NS}.
\end{equation}
We should state, that \eqref{NS} for a given sufficiently smooth right hand side togheter with sufficiently smooth initial values and boundary conditions admits a unique solution. We refer to the book of Temam \cite{T01} for more details. Whenever we want to emphasize the dependence of the solution on the control $u$ we will write $y=y(u)$.



\section{Dynamic Programming equation}

We illustrate the dynamic programming approach for abstract optimal control problems of the form
\begin{equation}\label{opt_con}
\min_{u\,\in\,\mathcal{U}}J_x(u):=\int_0^{\infty} L(y(t),u(t))\, e^{-\lambda t}\,dt\text{ subject to  } \dot y(t)=f(y(t),u(t)),\; y(0)=x,
\end{equation}
with system dynamics in  $\mathbb{R}^n.$  We assume $\lambda> 0$, and $L(\cdot,\cdot)$ and $f(\cdot,\cdot)$ to be Lipschitz-continuous, bounded functions. Then, it is clear that the optimal control problem \eqref{ocp} fits into the more abstract setting \eqref{opt_con}.\\
\noindent
In this setting, a standard solution tool is the application of the dynamic programming principle, which leads to a characterization of the value function $v(x):=\inf\limits_{u\in\mathcal{U}} J_x(u)$ as a viscosity solution of the Hamilton-Jacobi-Bellman equation (HJB)

\begin{equation}
\lambda v(x) - \inf_{u\in U}\{Dv\cdot f(x,u)+L(x,u)\}=0\,.\label{HJB}
\end{equation}

\noindent To approximate equation \eqref{HJB}, we construct a fully-discrete semi-Lagrangian scheme which is based on a discretization of the system dynamics with time step $h$, and a finite element discretization of the state space with mesh parameter $k$, leading to a fully discrete approximation $V_{h,k}(x)$ of the value function $v$ satisfying
\begin{equation}
V_{h,k}(x_i)=\min_{u\,\in U}\{(1-\lambda h)I_1[V_{h,k}](x_i+hf(x_i,u))+L(x_i,u)\}\,,\label{HJBh}
\end{equation}
for every element $x_i$ of the discretized spatial domain. In general, the arrival point $x_i+hf(x_i,u)$ is not a node of the state space grid, and therefore the value of $V_{h,k}$ at this point is approximated by means of a first-order interpolant of the data, denoted by $I_1[V_{h,k}]$ (we refer the reader to \cite[Appendix A]{BCD:97} for more details).

The goal is to find a feedback control law of the form $u(t)=\Phi(y(t),t)$ which steers the system to the desired trajectory. $\Phi$ is called {\em feedback map}. The computation of feedback maps is almost built in and comes straightforward from the knowledge of the value function. In fact;
$$\Phi(y_x(t))=u^*(t)=\arg\min_{u\in U} \left\{L(x,u)+\nabla v(x)^Tf(x,u)\right\}.$$
The characterization of the value function is valid for all classical problems in any dimension and its approximation is based on a-priori error estimates in $L^\infty$. \\
The request to solve an HJB in high dimensions comes up naturally whenever we want to control evolutive PDEs. However, a direct discretization, in many practically relevant situations, is impossible since the system of ODEs associated to a semi-discretization in time would have the dimension equal to the space dimension where one should solve the HJB equation.
Fortunately, at the discrete level, the POD (\cite{Sir87,V07})) method allows us to obtain low-dimensional reduced models even for complex dynamics, and, thus, presents an oppurtunity to circunmvent the curse of dimensionality in the numerical solution of the HJB equation.

\section{POD-Model Reduction for the controlled problem}

The Reduced Order Modelling (ROM) approach to optimal control problems is based on projecting the nonlinear dynamics onto a low dimensional manifold utilizing projectors that contain informations of the expected controlled flow. A common approach here is based on the snapshot form of POD proposed by Sirovich in \cite{Sir87}, which in the present situation works as follows.
We compute the snapshots set $y_1,\dots,y_n$ of the flow corresponding to different time instances $t_1,\ldots,t_n$ and define the POD ansatz of order $\ell$ for the state $y$ by

\begin{equation}\label{pod_ans}
y^\ell=\bar{y}+\sum_{i=1}^\ell w_i\psi_i,
\end{equation}


\noindent
where $\bar{y}=\frac{1}{n}\sum_{i=1}^n y_i$ denotes the mean flow and the basis functions $\{\psi_i\}_{i=1}^\ell$ are obtained from the singular value decomposition of the snapshot matrix $Y=[y_1-\bar{y},\ldots,y_n-\bar{y}],$ i.e. $Y=\Psi\Sigma V$, and the first $\ell$ columns of $\Psi$ form the POD basis functions of rank $\ell$.  Here the SVD is based on the Euclidean inner product. This is reasonable in our situation, since the numerical computations performed in our numerical example for the driven cavity problem are based on a uniform staggered grid. The snapshots are computed on the basis of a stable finite difference discretization of \eqref{NS} which leads to a semi-discretet system of ODEs of the form
\begin{equation}\label{NS_dis}
\dot{y}+\nu A y + Cp=\eta(y)+B u,\quad y(0)=y_0.
\end{equation}

\noindent
The reduced optimal control problem is obtained through replacing \eqref{NS_dis} by a dynamical system obtained from a Galerkin approximation with basis functions $\{\psi_i\}_{i=1}^\ell$ and ansatz \eqref{pod_ans} for the state.

\noindent
This leads to a $\ell-$dimensional system for the unknown coefficients $\{w_i\}_{i=1}^\ell,$ namely
\begin{equation}\label{NS_pod}
M^\ell\dot{w}+\nu A^\ell w=\eta(w)+B^\ell u\quad w(0)=w_0.
\end{equation}
Here the entries of the mass $M^\ell$ and the stiffness $A^\ell$ are given by $\langle\psi_j,\psi_i\rangle$ and $\langle\psi_j,A\psi_i\rangle$,  respectivelly. The reduced shape function is obtained by $(B^\ell)_i=\langle B,\psi_i\rangle.$ The coefficients of the initial condition $y^\ell(0)\in\R^\ell$ are determined by $w_i(0)=(w_0)_i=\langle y_0-\bar{y},\psi_i\rangle, \;\;1\leq i\leq \ell,$ and the solution of the reduced dynamical problem is denoted by $w(s)\in\R^\ell.$ Note that for the reduction of the nonlinear term $\eta(w)$ we use the {\em Discrete Empirical Interpolation Method} (DEIM, see\cite{CS09}).  The pressure does not appear in the reduced problem \eqref{NS_pod} since the snapshots are divergence-free.
Then, the POD-Galerkin approximation leads to  the optimization problem
\begin{equation}\label{KPl}
\inf J^\ell_{w_0} (u),\\
\end{equation}
where $u\in\mathcal{U}$, $w$ solves \eqref{NS_pod} and the cost functional is defined by
$$J^\ell_{w_0}(u)=\int_0^\infty L(w(s),u(s),s)e^{-\lambda s}\; ds.$$
The value function $v^\ell$, defined for the initial state $w_0\in\R^\ell$  is given by
$$v^\ell(w_0)=\inf_{u\in\mathcal{U}} J^\ell_{w_0}(u),$$ 
and $w$ solves (\ref{NS_pod}) with the control $u$ and initial condition $w_0.$ HJB equations are defined in $\R^n,$ but we need to restrict our numerically domain to a bounded subset of $\R^n$. We refer the interested reader to \cite{AF12} for a detailed description.\\

\noindent
\section{Numerical Tests}
\begin{figure}[htbp]
\centering
\includegraphics[scale=.24]{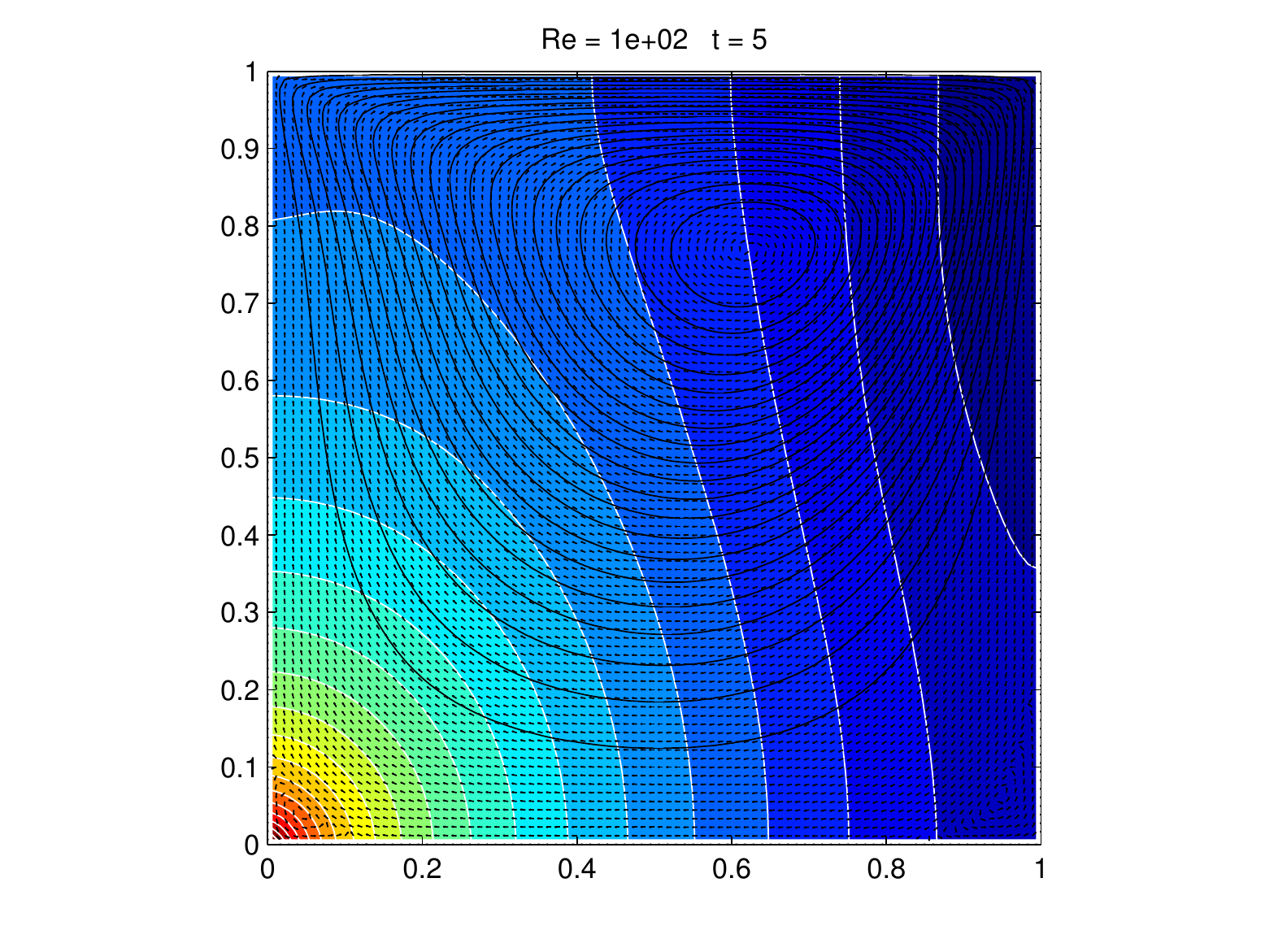}\hfill\includegraphics[scale=.24]{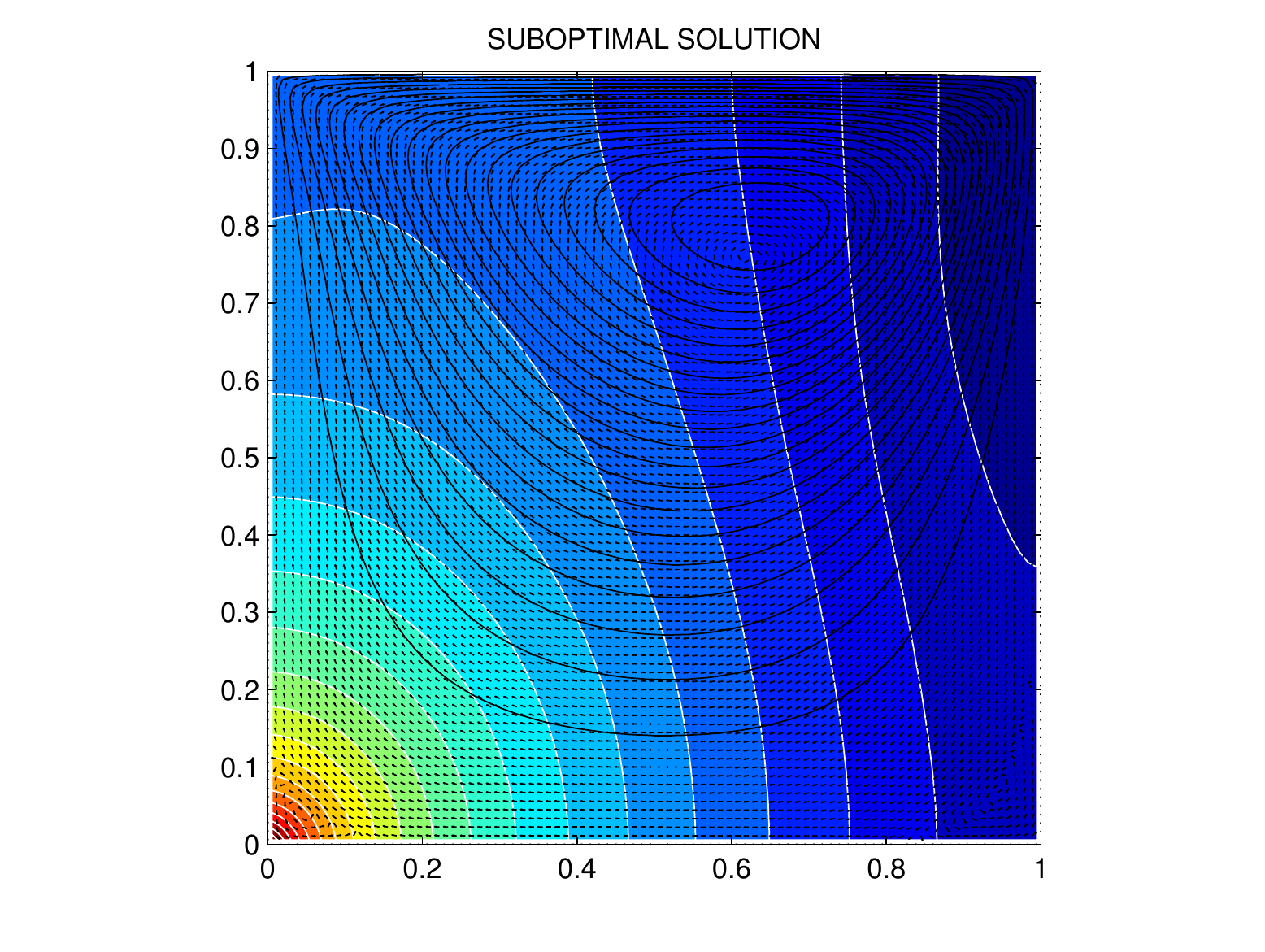}\hfill\includegraphics[scale=.24]{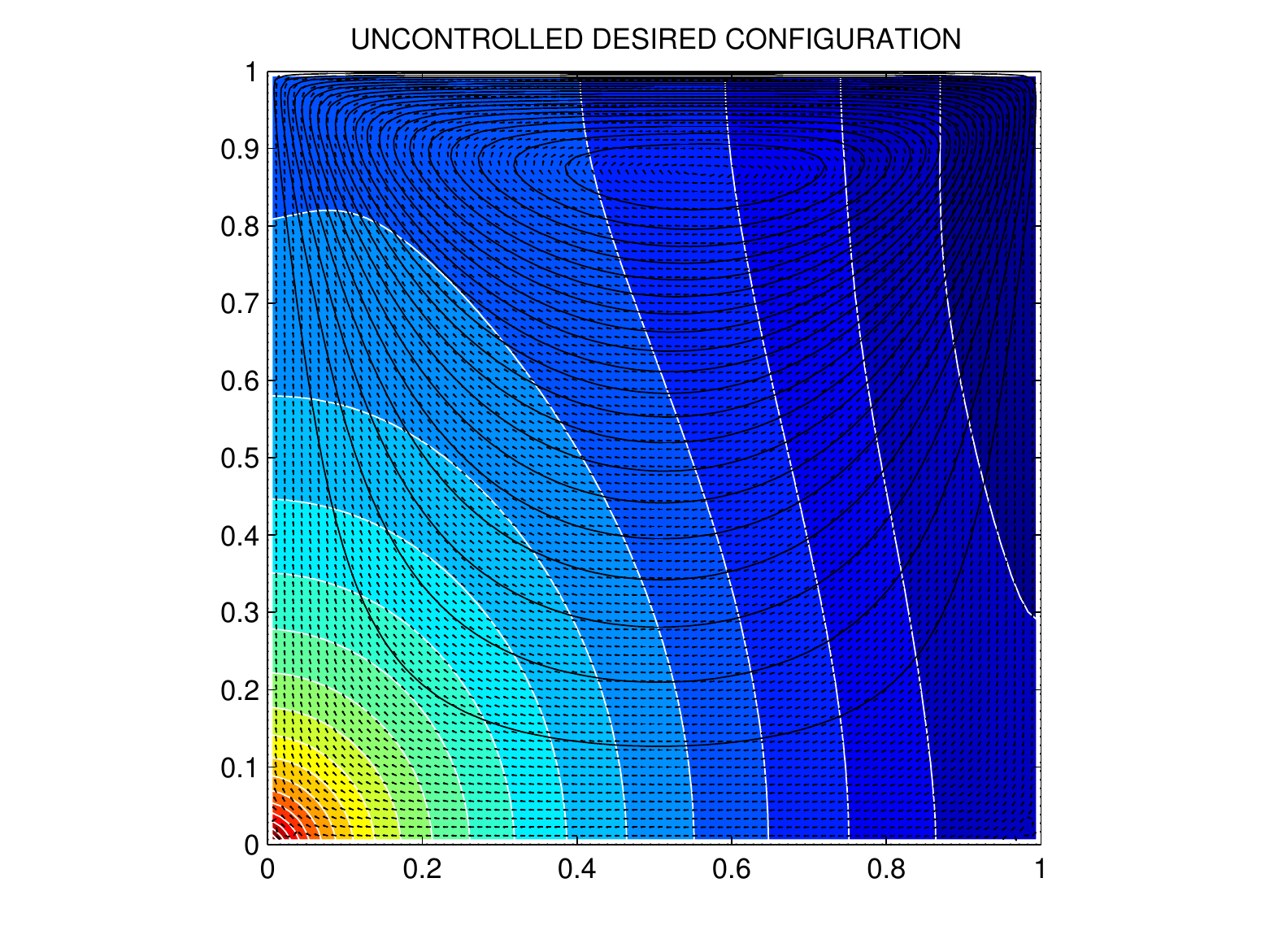}

\caption{Mean flow NS (left) - controlled configuration at time $t=0.5$ (middle) - uncontrolled configuration at time $t=0.5$ (right).}
\label{fig:des}
\end{figure}

In this section we consider as numerical example the control of the flow in the lid-driven cavity. 
In \eqref{ocp} we set: $\Omega=(0,1)\times(0,1), y_0\equiv0,\nu=0.01,\alpha=0.01, \lambda=1,U=\{-1,0,1\}, y_b=(1,0)$ on the top boundary and $y_b=(0,0)$ on the remaing boundary segments. In \eqref{HJBh} we take $k=0.2,$ $h=0.04$ whereas the optimal trajectory is obtained with a time stepsize of $0.01$.\\ 
The control gain of the suboptimal control problem, with the ansatz \eqref{pod_ans}, consists of steering the coefficients $w$ to the origin. For the purpose of this test, we take only 3 POD and 6 POD-DEIM basis functions. In our numerical computations the reduction of the nonlinearity with DEIM already yields a considerable computational speedup. Further investigations on the performance of DEIM in relation to the discretization parameters are provided in a subsequent paper. The snapshots are computed with a finite difference scheme from the uncontrolled problem ($u\equiv0$) in \eqref{NS} where we use the Matlab code provided in \cite{S08}. \\
\noindent
In Figure \ref{fig:des} we show the configuration of the flow. On the left we show the mean flow, which is the desired state, in the middle the controlled flow is shown, and on the right the uncontrolled flow is shown. As shape function we use the steady state solution of the Navier-Stokes system.

\noindent
We can see that at time $t=0.5$ the suboptimal solution already well approximates the desired state, as confirmed in Table \ref{table1}, where the $L^\infty-$error of $y^\ell-\bar{y}$ at $t=0.5$ and $t=4$ is reported for this shape function. When the time is increasing the solution itself tends to stabilize close to the mean flow, but still the suboptimal solution has a smaller error with respect to the uncontrolled problem. Note that the performance of our method depends on the choice of the shape functions. 
In Table \ref{table2} we display the results obtained with the steady state solution of the Stokes equation as shape function. As expected, the approach works better if we can use the steady state of the Navier-Stokes equation as shape function.

\begin{table}[htbp]
\begin{center}
\begin{tabular}{ccc}
\toprule

  & $t=0.5$ & $t=4$ \\
\midrule

$\|y^\ell(x,t,u^\ell)-\bar{y}\|_\infty$ & 0.007 & 0.006   \\

$\|y(x,t;0)-\bar{y}\|_\infty$ & 0.283 & 0.048\\
\bottomrule

\end{tabular}
\end{center}
\caption{$L^\infty$ error at time $t=0.5$ and $t=4.$ $\bar{y}$ is the desidered state, $y^\ell(x,t;u^\ell)$ is the suboptimal solution, and $y(x,t;0)$ denotes the uncontrolled solution. The shape function is chosen as the steady state solution of the Navier-Stokes equations.}
\label{table1}
\end{table}

\begin{table}[htbp]
\begin{center}
\begin{tabular}{ccc}
\toprule

  & $t=0.5$ & $t=4$ \\
\midrule

$\|y^\ell(x,t,u^\ell)-\bar{y}\|_\infty$ & 0.081 & 0.022   \\

$\|y(x,t;0)-\bar{y}\|_\infty$ & 0.283 & 0.048\\
\bottomrule

\end{tabular}
\end{center}
\caption{$L^\infty$ error at time $t=0.5$ and $t=4.$ $\bar{y}$ is the desidered state, $y^\ell(x,t;u^\ell)$ is the suboptimal solution, and $y(x,t;0)$ denotes the uncontrolled solution. The shape function is chosen as the steady state solution of the Stokes equations.}

\label{table2}
\end{table}

\begin{figure}[htbp]
\centering
\includegraphics[scale=.35]{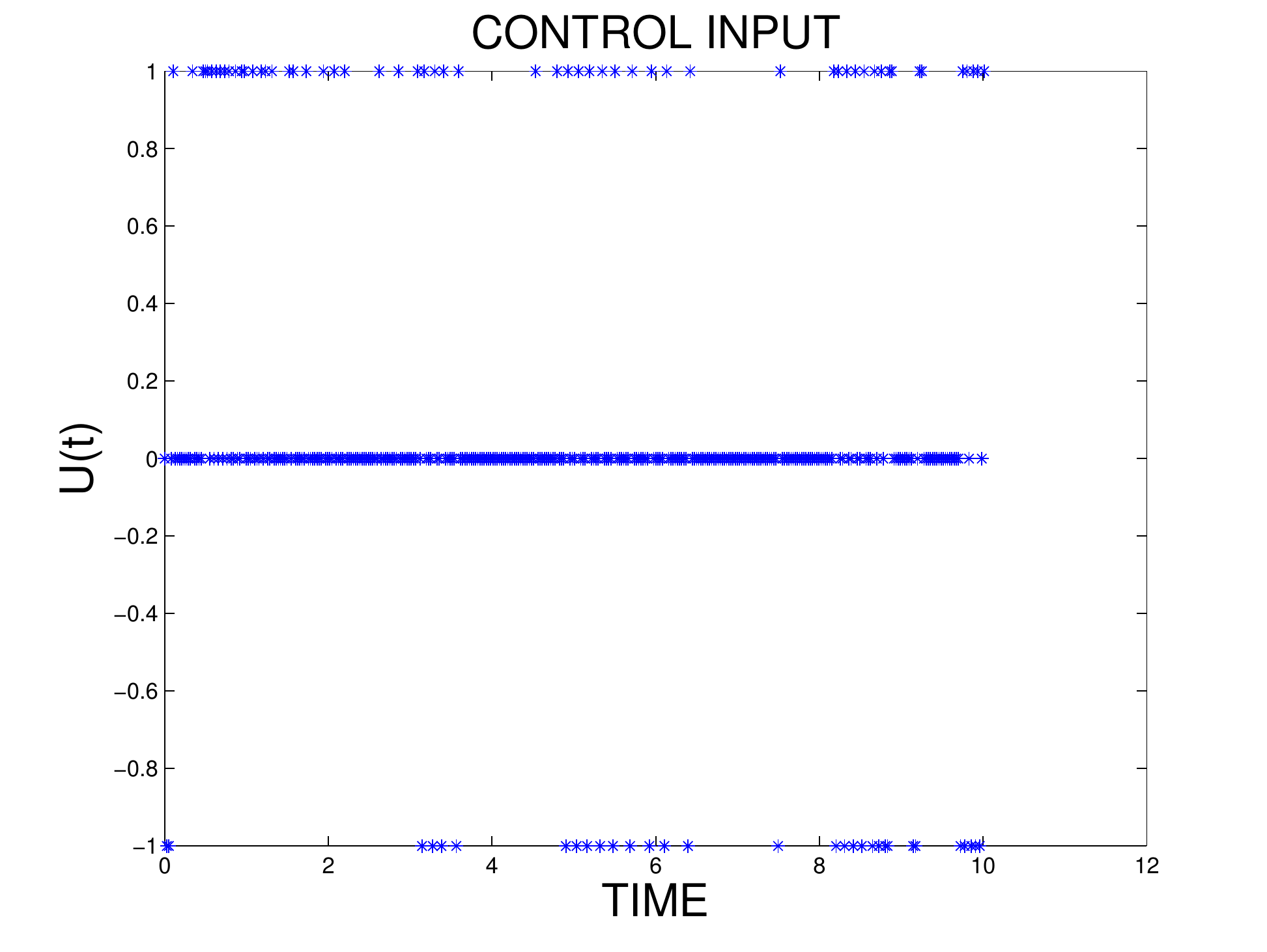}

\caption{Control input with 3 constant controls $\{-1, 0, 1\}$ and one shape function chosen as the steady state solution of the Navier-Stokes equation. }
\label{fig:vel}
\end{figure}

\noindent
In Figure \ref{fig:vel} we present the control input. The behavior of the control is classical for feedback control, since the system tries to correct step by step the trajectories. The control space is only given by constant values  $\{-1, 0, 1\}.$

%
%
%
%

%
%

%
%

%
%
\bibliographystyle{spmpsci}

\begin{thebibliography}{1}
\bibitem{AF12} A. Alla, M. Falcone. {\em An adaptive POD approximation method for the control of advection-diffusion equations} International Series of Numerical Mathematics (Birkhauser, Basel, 2013)
\bibitem{AF13b} A. Alla, M. Falcone, {\em A Time-Adaptive POD Method for Optimal Control Problems,} to appear in the Proceedings of the 1st IFAC Workshop on Control of Systems Modeled by Partial Differential Equations,

\bibitem{AK01} J.A. Atwell, B.B. King, {\em Proper Orthogonal Decomposition for Reduced Basis Feedback Controllers for Parabolic Equations}, Matematical and computer modelling. {\bf 33} (2001), 1-19.

\bibitem{BCD:97} M. Bardi, I. Capuzzo Dolcetta. { \em Optimal control and viscosity solutions
 of Hamilton-Jacobi-Bellman equations}. Birkhauser, Basel, 1997.
to appear on SIAM J. Sci. Comp. 
\bibitem{CS09} S. Chaturantabut, D.C. Sorensen, {\em Discrete Empirical Interpolation for NonLinear Model Reduction,} SIAM J. of Scientific Computing, {\bf 32} (2010), 2737-2764.
\bibitem{FF12} M. Falcone, R. Ferretti. Semi-Lagrangian Approximation Schemes for Linear and 
Hamilton-Jacobi Equations, SIAM, 2013.
\bibitem{KV99} K. Kunisch, S. Volkwein. \emph{Control of Burgers' Equation by a Reduced Order Approach using Proper Orthogonal Decomposition}. Journal of Optimization Theory and Applications, 102 (1999), 345- 371.
\bibitem{KV01} K. Kunisch, S. Volkwein. \emph{Galerkin proper orthogonal decomposition 
methods for parabolic problems} Numer. Math. 90 (2001), 117-148.
\bibitem{KVX04} K. Kunisch, S. Volkwein, L. Xie. \emph{HJB-POD Based Feedback Design for the Optimal Control of Evolution Problems}. SIAM J. on Applied Dynamical Systems, 4 (2004), 701-722.
\bibitem{KX05} K. Kunisch, L. Xie. {\em POD-Based Feedback Control of Burgers Equation by Solving the Evolutionary HJB Equation}, Computers and Mathematics with Applications. {\bf 49} (2005), 1113-1126.
\bibitem{S08} B. Seibold, {\em A compact and fast Matlab code solbing the incompressible Navier-Stokes equations on rectangular domains}, 2008\\ {\tt http://math.mit.edu/cse/codes/mit18086-navierstokes.pdf}
\bibitem{Sir87} L. Sirovich, {\em Turbulence and the dynamics of coherent structures. Parts I-II,}
Quarterly of Applied Mathematics, {\bf XVL} (1987), 561-590.
\bibitem{T01} R. Temam, {\em Navier-Stokes Equations: Theory and Numerical Analysis}, American Mathematical Society 2001.
\bibitem{V07} S. Volkwein, {\em Model Reduction using Proper Orthogonal Decomposition}, 2011  {\tt www.math.uni-konstanz.de/numerik/personen/\\volkwein/index.php}
 \end{thebibliography}



\end{document}